\documentclass[a4paper,11pt,leqno]{amsart}

\usepackage{amsmath,amsfonts,amssymb}
\usepackage{pifont,graphicx,epsfig}
\usepackage[bf]{caption2}
\usepackage{epstopdf}
\usepackage{CJK}
\usepackage{scrextend}
\def\C{\hbox{\font\dubl=msbm10 scaled 1100 {\dubl C}}}
\def\Z{\hbox{\font\dubl=msbm10 scaled 1100 {\dubl Z}}}
\def\R{\hbox{\font\dubl=msbm10 scaled 1100 {\dubl R}}}
\def\N{\hbox{\font\dubl=msbm10 scaled 1100 {\dubl N}}}

\def\Re{{\rm{Re}}\,}

\newtheorem{Theorem}{Theorem}
\newtheorem{Corollary}[Theorem]{Corollary}
\newtheorem{Lemma}[Theorem]{Lemma}

\title[The Riemann Zeta-Function on Discrete Sets]{The Values of the Riemann Zeta-Function on Discrete Sets}

\author[J. Lee, A. Sourmelidis, J. Steuding,  A.I. Suriajaya]{Junghun Lee, Athanasios Sourmelidis,\\ J\"orn Steuding, Ade Irma Suriajaya}

\date{1st of November 2017}

\dedicatory{Dedicated to Professor Kohji Matsumoto at the occasion of his 60th Birthday}
\begin{CJK}{UTF8}{min}
\begin{document}

\maketitle

\begin{abstract}
We study the values taken by the Riemann zeta-function $\zeta$ on discrete sets. We show that infinite vertical arithmetic progressions are uniquely determined by the values of $\zeta$ taken on this set. Moreover, we prove a joint discrete universality theorem for $\zeta$ with respect to certain permutations of the set of positive integers. Finally, we study a generalization of the classical denseness theorems for $\zeta$.     
\end{abstract}

{\small \noindent {\sc Keywords:} Riemann zeta-function, value-distribution, universality}

\bigskip

\section{Introduction and Statement of the Main Results}

In view of the extensive research for about one century and {\bf 60 years} it might be surprising how little we know about the distribution of values of the Riemann zeta-function, which is, for $\Re s>1$, defined by 
$$
\zeta(s)=\sum_{n\geq 1}n^{-s}=\prod_p(1-p^{-s})^{-1},
$$
and by analytic continuation elsewhere except for a simple pole at $s=1$. The famous yet unsolved Riemann hypothesis claims that all nontrivial (non-real) zeros of $\zeta$ lie on the so-called critical line ${1\over 2}+i\R$ (which is equivalent to the non-vanishing of $\zeta(s)$ in the half-plane $\Re s>{1\over 2}$). This conjecture is one of the six unsolved millennium problems and dates back to the fundamental work \cite{rie} of Bernhard Riemann who was the first to study $\zeta$ in 1859 as a function of a complex variable and to contribute significant ideas for further investigations undertaken by hundreds of mathematicians eversince. It was first shown by Godfrey Hardy \cite{hardy} that there are infinitely many zeros on the critical line and the present best quantitative result in this direction is due to Feng \cite{feng} who proved, building on important work by Selberg, Levinson and others, that more than $41.28$ percent of the zeros are on the critical line. Except the value $0$ so far no other {\it explicit} complex number is known to be taken as a value by $\zeta({1\over 2}+it)$ for some real number $t$. Another challenging open question is whether the zeta-values taken on the critical line are dense in the complex plane or not. 

To the right of the critical line the situation is sligthly better. By the work of Bohr \& Courant \cite{bc} it is known that the sets $\zeta(\sigma+i\R)$ (defined as the set of values $\zeta(\sigma+it)$ for $t\in\R$) are dense in $\C$ for any fixed $\sigma\in({1\over 2},1]$, whereas for $\sigma>1$ there is no densenes (since the absolute convergence of the $\zeta$ defining series implies $\vert \zeta(s)\vert \leq \zeta(\sigma)$). Concerning the distribution of zeros it is a classical result that $\zeta(s)$ does not vanish for $\Re s\geq 1$ (which already implies the prime number theorem without remainder term). The present best zero-free region had been obtained by Vinogradov \cite{vino} and (independently) Korobov \cite{koro} and it is some funnel-shaped domain with real part asymptotically tending to $1$ from the left, too complicated to be reproduced here. Moreover, there cannot be too many exceptional zeros off the critical line as already follows from the classical density theorem due to Bohr \& Landau \cite{bl} who showed that the proportion of nontrivial zeros $\rho=\beta+i\gamma$ satisfying $\beta>\sigma$ for some fixed $\sigma>{1\over 2}$ and $0<\gamma\leq T$ (counting multiplicities) is vanishing as $T$ tends to infinity. Since nontrivial zeros come in pairs, symmetrically distributed with respect to the critical line (thanks to the functional equation), consequently, almost all zeros of the zeta-function are clustered around the critical line.  

Interestingly, the same clustering holds for $a$-points as well. We recall the definition of $a$-points. Given a complex number $a$, the roots of the equation 
$$
\zeta(s)-a=0
$$
are called $a$-points. For every $a$, they come in numbers similar to the zeros and some aspects of their distribution are pretty much the same. Landau \cite{landau} proved that almost all $a$-points are clustered around the critical line provided the Riemann hypothesis is true. The latter assumption was removed by Levinson \cite{levia} who showed (in a quantitative way) that almost all $a$-points are clustered around the critical line.

Maybe the most remarkable result about the value-distribution of the zeta-function, however, is Voronin's celebrated universality theorem which, roughly speaking, states that any non-vanishing analytic function $f$, defined on a closed disk $K$ inside the open right half of the critical strip, can be approximated by certain vertical shifts of the zeta-function; more precisely, for every positive $\epsilon$, there exists a real $\tau$ such that 
$$
\max_{s\in K}\vert \zeta(s+i\tau)-f(s)\vert <\epsilon;
$$
it is implicit in Voronin's proof \cite{voronin} that the set of shifts $\tau$ has even positive lower density. This powerful concept of approximation may be used to derive the denseness theorems mentioned above (by constant target functions $f$) whenever ${1\over 2}<\Re s<1$ and quite a few more interesting results. In the meanwhile several extensions and generalizations of Voronin's universality theorem have been achieved, in particular a discrete version (with shifts from an arithmetic progression) by Reich \cite{reich}, and in place of disks $K$ may also denote any compact subset of the right open half of the critical strip with connected complement. This and some further results concerning the values of the zeta-function will enter the stage in the following sections. For more information about the value-distribution of the zeta-function and its relatives, however, we refer to Kohji Matsumoto's excellent survey articles \cite{kohji1,kohji2}. 
\smallskip

In this note we shall investigate the values taken by the zeta-function on discrete subsets of the complex plane. As everyone knows, a subset $M$ of the complex plane is discrete if for every element $m\in M$ there exists a neighbourhood which does not contain any element from $M$ except $m$. Obviously, a discrete set is not dense. Typical examples of discrete sets are (finite or infinite) arithmetic progressions. In view of the aforementioned results we may restrict to vertical arithmetic progressions. Our first result is located in the half-plane of absolute convergence of the $\zeta$ defining Dirichlet series.

\begin{Theorem}\label{absconv}
Let $t_1,t_2$ be arbitrary real numbers and $\delta_1,\delta_2$ be arbitrary positive real numbers. Assume 
$$
\zeta(s+i(t_1+\delta_1n))=\zeta(s+i(t_2+\delta_2\sigma(n)))\qquad\mbox{for}\quad n=1,2,\ldots,
$$
where $\sigma\,:\,\N\to\N$ is bijective. Then, we have $t_1=t_2, \delta_1=\delta_2, \sigma={\rm{id}}$ or ${\rm{Re}}\,s\leq b$, where $b$ is a constant, explicitly given by (\ref{neq}) and depending only on $t_j,\delta_j$ for $j=1,2$ and $\sigma$.
\end{Theorem}

\noindent We may interpret this result as follows: an infinite vertical arithmetical progression is characterized by the values of the zeta-function taken at the values generated by this arithmetic progression. For short: for sufficiently large $\Re s$,
$$
\zeta(s+i(t_1+\delta_1\N))=\zeta(s+i(t_2+\delta_2\sigma(\N)))\quad\Rightarrow\quad t_1=t_2, \delta_1=\delta_2, \sigma={\rm{id}}.
$$

In the next result we consider the situation in the right half of the critical strip (in a sligthly more general context). For a real number $\alpha$ the Beatty sequence $\mathcal{B}_\alpha$ is defined as the sequence (or sometimes set) of numbers $\lfloor n\alpha\rfloor$ for $n\in\mathbb{N}$, where $\lfloor x\rfloor$ denotes the largest integer which is less than or equal to $x$. Of particular interest is the case of irrational $\alpha$. The well-known Beatty's (or in some literature Rayleigh's) theorem states that, for irrational $\alpha>1$,  
\begin{equation}\label{ray}
\mathcal{B}_\alpha\cap\mathcal{B}_{\alpha'}=\emptyset\hspace{0.2cm}
\text{and}\hspace{0.2cm}\mathcal{B}_\alpha\cup\mathcal{B}_{\alpha'}=\mathbb{N},
\end{equation}
where $\alpha'$ is determined by 
$$
\dfrac{1}{\alpha}+\dfrac{1}{\alpha'}=1. 
$$
We shall use this result to generate infinitely many permutations of the set of positive integers as follows. For every irrational $\alpha>1$ define a bijection $\sigma_\alpha:\mathbb{N}\to\mathbb{N}$ by
\begin{equation}\label{permutation}
n\mapsto \sigma_\alpha(n)=\left\{\begin{array}{ll}\lfloor m\alpha'\rfloor,&\text{if }n=\lfloor m\alpha\rfloor\text{ for some }m\in\mathbb{N},\\
\lfloor m\alpha\rfloor,&\text{if }n=\lfloor m\alpha'\rfloor\text{ for some }m\in\mathbb{N}.
\end{array}\right.
\end{equation}
We consider the values taken by the zeta-function on discrete sets formed with respect to such permutations $\sigma$ of $\mathbb{N}$. For this purpose we define, for $\delta_1$ and $\delta_2$ arbitrary positive real numbers, the sets
$$
\mathcal{A}_i=\left\{\delta_i\dfrac{\log q}{2\pi}:q\in\mathbb{Q}_+\right\},\hspace*{0.15cm}i=1,2,
$$
and 
$$
\mathcal{L}(\delta_1,\delta_2)=\Bigg\lbrace \alpha\in\mathbb{R}\hspace*{-0.05cm}\,:\,\hspace*{-0.05cm}\forall( \theta_1, \theta_2)\hspace*{-0.05cm}\in\hspace*{-0.05cm}\mathcal{A}\ \hspace*{-0.05cm} :\hspace*{-0.05cm}
1,\alpha,\alpha',\alpha\theta_1+\alpha'\theta_2\hspace*{0.15cm}\text{are lin. ind. over }\mathbb{Q}\Bigg\rbrace,
$$
where $\mathcal{A}:=\mathcal{A}_1\times\mathcal{A}_2\setminus\lbrace(0,0)\rbrace.$
Finally, we denote by $\mathcal{D}$ the right half of the critical strip, i.e., $\mathcal{D}=\{s\in\C\,:\,{\rm{Re}}\,s\in({1\over 2},1)\}$. Then  

\begin{Theorem}\label{thana}Let $t_1,t_2$ be arbitrary real numbers and $\delta_1,\delta_2$ be arbitrary positive real numbers.
Let also $\alpha\in\mathcal{L}(\delta_1,\delta_2)\cap(1,+\infty),$ $K$ be a compact subset of $\mathcal{D}$ with connected complement and $f,g$ continuous non-vanishing functions on $K,$ which are analytic in the interior of $K.$ Then, for every $\varepsilon>0,$ 
$$
\liminf\limits_{N\to\infty}\dfrac{1}{N}\#\left\{1\leq n\leq N:\begin{array}{ll}
\max\limits_{s\in K}|\zeta(s+i(t_1+\delta_1\lfloor n\alpha\rfloor))-f(s)|<\varepsilon\\
\max\limits_{s\in K}|\zeta(s+i(t_2+\delta_2\lfloor n\alpha'\rfloor))-g(s)|<\varepsilon
\end{array}\right\}>0.\hspace*{1cm}
$$
\end{Theorem}

\noindent This may be considered as a joint discrete universality theorem with shifts taken with respect to the dissection (\ref{ray}) of $\mathbb{N}$ generated by the Beatty sequences of $\alpha$ and $\alpha'$. As an immediate consequence, we have   

\begin{Corollary}\label{sis}Let $t_1,t_2$ be arbitrary real numbers and $\delta_1,\delta_2$ be arbitrary positive real numbers.
Let also $\alpha\in\mathcal{L}(\delta_1,\delta_2)\cap(1,+\infty),$ $\sigma_\alpha$ be the permutation defined by (\ref{permutation}), $K$ be a compact subset of $\mathcal{D}$ with connected complement and $f,g$ continuous non-vanishing functions on $K,$ which are analytic in the interior of $K.$ Then, for every $\varepsilon>0,$
$$\liminf\limits_{N\to\infty}\dfrac{1}{N}\#\left\{1\leq n\leq N:\begin{array}{ll}
\max\limits_{s\in K}|\zeta(s+i(t_1+\delta_1n))-f(s)|<\varepsilon\\
\max\limits_{s\in K}|\zeta(s+i(t_2+\delta_2\sigma_\alpha(n)))-g(s)|<\varepsilon
\end{array}\right\}>0.$$
\end{Corollary}

\noindent Thus, using different target functions, we have, for ${\rm{Re}}\,s\in({1\over 2},1)$,
$$
\zeta(s+i\N)\neq \zeta(s+i\sigma_\alpha(\N))
$$
with $\sigma_\alpha$ as in the corollary; one could extend this result to arbitrary arithmetic progressions (of the form as in Theorem \ref{absconv}), however, one would have to consider quite a few cases (with respect to the common differences $\delta_j$) so we restrict on this simple form for the sake of simplicity.
\par

Finally, we consider the distribution of values taken on arithmetic progressions with respect to rather general sets. This extends the question about the denseness of values of $\zeta$. For this purpose define
$$
\mathbb{M} = \{  \mathsf{M} \subset \hat{\mathbb{C}} \mid \mathsf{M}^{\circ} \neq \emptyset \},
$$
where $A^{\circ}$ denotes the interior of a set $A$ in $\hat{\mathbb{C}}:=\mathbb{C}\cup\{\infty\}$. Moreover, let $\mathbb{M}_{\infty} = \{  \mathsf{M} \in \mathbb{M} \mid \infty \in \mathsf{M}^{\circ} \}$. Then

\begin{Theorem}\label{AZTHM1}
Let $h$ be a positive real number (we write $h \in \mathbb{R}_{>0}$), $l \in \mathbb{N}$, and $\mathsf{M} \in \mathbb{M}$. If $s \in \mathbb{C}$ satisfies $\operatorname{Re}(s)\in (1/2, 1)$, then there exists an infinite subset $\mathsf{N} \subset \mathbb{N}$ such that
$$
\zeta(s + i h (n + k - 1)) \in \mathsf{M}
$$
for any $n \in \mathsf{N}$ and $k \in \{ 1, 2, \cdots, l \}$.
\end{Theorem}

\noindent Notice that this theorem also implies that any bounded set does not contain any set of values of $\zeta$ on a vertical arithmetic progression of infinite length. We can see this by considering the complement of the set $\mathsf{M}$ in $\hat{\mathbb{C}}$, that is $\hat{\mathbb{C}}\backslash\mathsf{M}$, in Theorem \ref{AZTHM1} when $\infty$ is in $\mathsf{M}^{\circ}$. Then it follows that $\zeta(s + i h (n + k - 1)) \in \mathsf{M}$, therefore $\zeta(s + i h (n + k - 1))$ is not in the bounded set $\hat{\mathbb{C}}\backslash\mathsf{M}$.
\par

We can also flip the result of Theorem \ref{AZTHM1} to the left half of the critical strip by using the functional equation when $\infty$ is an interior point of the subset of $\hat{\mathbb{C}}$.

\begin{Corollary}\label{AZTHM2}
Let $h \in \mathbb{R}_{>0}, l \in \mathbb{N}$, and $\mathsf{M} \in \mathbb{M}_{\infty}$.
If $s \in \mathbb{C}$ satisfies $\operatorname{Re}(s) \in (0, 1/2)$, then there exists an infinite subset $\mathsf{N} \subset \mathbb{N}$ such that
$$
\zeta(s + i h (n + k - 1)) \in \mathsf{M}
$$
for any $n \in \mathsf{N}$ and $k \in \{ 1, 2, \cdots, l \}$.
\end{Corollary}

\noindent Putnam \cite{Putn54a,Putn54b} showed that the sequence of consecutive nontrivial zeros of $\zeta$ on the critical line does not contain any arithmetic progression of infinite length.
His result can be interpreted using our notation as considering $l = 1$, $s = 1/2$, and $0 \notin \mathsf{M}$.
Using fractal zeta-functions, Lapidus and Frankenhuijsen \cite[Chapter $11$]{LF13} gave a new proof of Putnam's result and extended it to a large class of zeta functions and $L$-series.
Our approach is based on the universality of the Riemann zeta function and is different from both of theirs.
However, we remark that our results do not include their result since $\operatorname{Re}(s)\neq1/2$ in Theorem \ref{AZTHM1} and Corollary \ref{AZTHM2}.
\smallskip

This paper is organized as follows. We begin in the half-plane of absolute convergence with the proof of Theorem \ref{absconv} in the following section, continue analytically to the right half of the critical strip and provide the proofs of Theorem \ref{thana}, its Corollary \ref{sis}, and Theorem \ref{AZTHM1}. Finally, we move further to the left, namely beyond the critical line, and prove Corollary \ref{AZTHM2} to the latter result by using the functional equation. 

\section{The Half-Plane of Absolute Convergence}

Our reasoning for Theorem \ref{absconv} applies in a more general context. Ordinary Dirichlet series converge (absolutely) in half-planes (if they converge). It was shown by Reich \cite{reich2} that the closure of the values of an ordinary Dirichlet series on a vertical line $\sigma+i\R$  inside the half-plane of absolute convergence is identical with the closure of its values taken on an arithmetic progression $\sigma+i\delta\N$ if, and only if, the numbers $1,{\delta\over 2\pi}\log p_1,\ldots, {\delta\over 2\pi}\log p_n$ are linearly independent over $\Z$, where $n$ is any positive integer and $p_n$ denotes the $n$th prime number in ascending order. Comparing two arithmetic progressions of this type, one should not expect the set of values to be identical without taking the closure. This is what we aim to investigate in the sequel.
\par

For $j=1,2$, let 
$$
L(s,f_j):=\sum_{m\geq 1}f_j(m)m^{-s}, 
$$
where the coefficients are bounded, say $\vert f_j(m)\vert \leq B$, whereby the Dirichlet series converge absolutely for $\Re s>1$. We shall consider the implications of the following equality of multisets,
\begin{equation}\label{multi1}
L(s+i(t_1+\delta_1\N),f_1)=L(s+i(t_2+\delta_2\N),f_2),
\end{equation}
where $t_j\in\R, 0<\delta_j\in\R$, and $L(s+i(t+\delta\N),f)$ is defined to be the multiset consisting of the elements $L(s+i(t+\delta n),f)$ for $n\in\N$. Hence, making use of a bijective mapping $\sigma\,:\,\N\to\N$, we may reformulate our set identity by the following equalities of values of the corresponding Dirichlet series:
\begin{equation}\label{multi2}
L(s+i(t_1+\delta_1n),f_1)=L(s+i(t_2+\delta_2\sigma(n)),f_2)\qquad\mbox{for}\quad n\in\N.
\end{equation}

For a positive integer $n$ we define
$$
\Phi_n(s)=L(s+i(t_1+\delta_1n),f_1)-L(s+i(t_2+\delta_2\sigma(n)),f_2)=\sum_{m\geq 1}\phi_n(m)m^{-s},
$$
where
\begin{align}\label{phi}
\phi_n(m)=f_1(m)m^{-i(t_1+\delta_1n)}-f_2(m)m^{-i(t_2+\delta_2\sigma(n))};
\end{align}
obviously, this Dirichlet series converges for $\Re s>1$. If $\Phi_n$ is not identically vanishing, by the uniqueness of Dirichlet series representation (see \cite{apostol}), there exists
$$
\mu:=\min\{m\in\N\,:\,\phi_n(m)\neq 0\}
$$
which depends on $n$ (although we do not indicate this here). In this case
$$
\Phi_n(s)=\phi_n(\mu)\mu^{-s}+\sum_{m\geq \mu+1}\phi_n(m)m^{-s}.
$$
Since
$$
\left\vert \sum_{m\geq \mu+1}\phi_n(m)m^{-s}\right\vert\leq 2B\int_\mu^\infty u^{-{\rm{Re}}\,s}{\rm{d}}\,u={2B\mu^{1-{\rm{Re}}\,s}\over {\rm{Re}}\,s-1},
$$
it follows that 
$$
\vert \Phi_n(s)\vert\geq \vert \phi_n(\mu)\vert\mu^{-{\rm{Re}}\,s} -{2B\mu^{1-{\rm{Re}}\,s}\over {\rm{Re}}\,s-1}.
$$
Hence, 
$$
L(s+i(t_1+\delta_1n),f_1)-L(s+i(t_2+\delta_2\sigma(n)),f_2)\neq 0
$$
for $\Re s>b_n:=1+{2B\mu\over \vert \phi_n(\mu)\vert}$. Now let 
\begin{equation}\label{neq}
\begin{array}{rr}
b:=\inf\{b_n\,:\,\phi_n(m)\neq 0\}.
\end{array}
\end{equation}
Then, for at least one $n$, we have 
$$
L(s+i(t_1+\delta_1n),f_1)\neq L(s+i(t_2+\delta_2\sigma(n)),f_2)\qquad\mbox{for}\quad \Re s>b,
$$
hence, (\ref{multi2}) and also (\ref{multi1}) cannot hold for sufficiently large $\Re s$ under the assumption of the existence of a positive integer $\mu$ such that $\phi_n(\mu)\neq 0$. Moreover, it follows that for any fixed $n$ the function $\Phi_n(s)$ is vanishing only on a discrete subset of the right half-plane ${\rm{Re}}\,s>1$, and for (\ref{multi1}) to hold all the functions $\Phi_n(s)$ for $n\in\N$ would need to have a common zero.   
\par

If $\vert f_1(m)\vert \neq\vert f_2(m)\vert$ for some $m\in\N$, then  $\phi_n(\mu)\neq 0$ for every $n\in\mathbb{N}$, as it can be seen from (\ref{phi}), and we may conclude (\ref{neq}). Hence, in the sequel we shall suppose that $\vert f_1(m)\vert =\vert f_2(m)\vert$ for all $m\in\N$.
\par

Suppose that $\phi_n(m)=0$ for all $m,n\in\N$. Then 
\begin{equation}\label{equi}
f_1(m)=f_2(m)m^{i(t_1+\delta_1n-(t_2+\delta_2\sigma(n)))}.
\end{equation}

If $f_1=f_2$, then $\exp(i\Delta_n\log m)=1$ with $\Delta_n:=t_1+\delta_1n-(t_2+\delta_2\sigma(n))$ whenever $f_1(m)\neq 0$. Assume that the latter condition holds for at least two coprime $m\neq 1$. Since the logarithms of the prime numbers are linearly independent over the rationals (as follows from unique prime factorization of the integers), it then follows that $\Delta_n=0$, or
$$
t_1+\delta_1n=t_2+\delta_2\sigma(n)\qquad\mbox{for all}\quad n\in\N.
$$
Thus,
\begin{eqnarray*}
\delta_1 &=&t_1+\delta_1(n+1)\,-\,(t_1+\delta_1n)\\
&=& t_2+\delta_2\sigma(n+1)-(t_2+\delta_2\sigma(n))\,=\delta_2(\sigma(n+1)-\sigma(n)).
\end{eqnarray*}
Hence, $\sigma(n+1)-\sigma(n)$ is constant which, since $\sigma$ is a permutation of $\N$, implies $\sigma={\rm{id}}$. Moreover, this yields $\delta_1=\delta_2$ and $t_1=t_2$. This shows that (\ref{multi1}) can hold for this trivial case only and so we have proved:

\begin{Theorem}\label{absconve}
Let $f$ be a bounded arithmetical function, supported on a set of positive integers having at least two distinct prime divisors, and let $t_1,t_2$ be arbitrary real numbers and $\delta_1,\delta_2$ be arbitrary positive real numbers. Assume 
$$
L(s+i(t_1+\delta_1n),f)=L(s+i(t_2+\delta_2\sigma(n)),f)\qquad\mbox{for}\quad n=1,2,\ldots,
$$
where $\sigma\,:\,\N\to\N$ is bijective. Then, we have $t_1=t_2, \delta_1=\delta_2, \sigma={\rm{id}}$ or ${\rm{Re}}\,s\leq b$, where $b$ is a constant, explicitly given by (\ref{neq}) and depending only on $f,t_j,\delta_j$, and $\sigma$.
\end{Theorem}

\noindent Of course, this implies Theorem \ref{absconv}. 
\par

Notice that the condition concerning the prime factors are necessary as follows from the following example
$$
(1-2^{-(s+2\pi i/\log 2)})^{-1}=(1-2^{-s})^{-1}=\sum_{j\geq 0}2^{-js}.
$$
Rewriting the Dirichlet series on the right hand side as $L(s,f)$, we observe that $f(n)\neq 0$ if, and only if, $n=2^j$ for some non-negative integer $j$, and the displayed formula of Theorem \ref{absconve} holds with $t_1={2\pi \over \log 2}\neq 0=t_2, \delta_1=\delta_2, \sigma={\rm{id}}$. 

Theorem \ref{absconve} implies the somewhat surprising result that if a value is taken twice by $L(s,f)$ on an arithmetic progression $s+i(t+\delta n)$, then the real part of $s$ has to be {\it small}. In fact, setting $\delta:=\delta_1=\delta_2, t:=t_1=t_2$ and assuming that there would exist distinct positive integers $n_1$ and $n_2$ satisfying 
$$
L(s+i(t+\delta n_1),f)=L(s+i(t+\delta n_2),f),
$$
would lead to a contradiction by defining $\sigma(n_1)=n_2, \sigma(n_2)=n_1$, and $\sigma(n)=n$ for all $n\neq n_1,n_2$. Of course, the curve $t\mapsto L(\sigma+it,f)$ may have intersections as $t$ ranges through $\R$, however, the condition (\ref{neq}) on the real part $\sigma$ excludes their appearance in arithmetic progression. It is natural to ask whether this result can be extended to the left. 
\par

As already mentioned in the introduction, it is unknown whether the curve $t\mapsto \zeta({1\over 2}+it)$ is dense in $\C$ or not. However, there are further interesting questions. Using a computer algebra package, one observes that this curve has infinitely many intersection points, in particular at the origin (corresponding to the zeta zeros on the critical line). One may ask whether for any fixed $c\in\zeta({1\over 2}+i\R)$ different from zero the number of real $t$ satisfying $c=\zeta({1\over 2}+it)$ is bounded or not; the same question for $t$ from an arithmetic progression might be easier to study. 
\par

If we consider the case of two distinct arithmetical functions $f_1\neq f_2$ with $f_1(m)=f_2(m)\neq 0$ for some $m\geq 2$, then it follows from (\ref{equi}) that 
$$
t_1+\delta_1 n\equiv t_2+\delta_2\sigma(n)\bmod\, 2\pi/\log m.
$$
If we assume that $t_1-t_2,\delta_1,\delta_2$ and ${2\pi \over \log m}$ are linearly independent over the rationals, we get a contradiction. For a general result one would have to discuss quite a few cases besides this generic case. We leave this to the interested reader and mention only that this theme is not unrelated to Bohr's notion of {\it equivalent} Dirichlet series \cite{bohr}.
\par

Concerning the discrete set of zeros of a function $\Phi_n(s)$ in case of the zeta-function, it follows from Nevanlinna's value distribution theory that 
$$
\Phi_n(s)=\zeta(s+i(t_1+\delta_1n))-\zeta(s+i(t_2+\delta_2\sigma(n)))
$$
being a linear combination of shifts of $\zeta(s)$ would have the same characteristic function as $\zeta(s)$ itself, which is $T(r,\zeta)={r\over \pi}\log r+O(r)$ (see \cite{nevanlinna,Ste07} for details). This implies that the number of zeros $\rho_n$ of $\Phi_n(s)$ satisfying $\vert\rho_n\vert\leq r$ is $O(r\log r)$ as $r\to\infty$. 

\section{The Right Half of the Critical Strip}

For every $m\in\mathbb{N},$ define the truncated Euler products $\zeta_m$ by the formula 
$$
\zeta_m(s)=\prod\limits_{i=1}^m(1-p_i^{-s})^{-1},
$$
where $s\in\left\{z\in\mathbb{C}:\mathrm{Re}z>0\right\}$ and $p_i$ denotes the $i$-th prime(in ascending order). As $m$ tends to infinity, these products do not converge in the critical strip $\mathcal{D}$ but they approximate $\zeta$ in the mean, that is
\begin{align}\label{m.v.a.c.}
\lim\limits_{m\to\infty}\limsup\limits_{T\to\infty}\dfrac{1}{T}\int_0^T|\zeta(s+i\tau)-\zeta_m(s+i\tau)|^2\mathrm{d}\tau=0
\end{align}
uniformly on compact subsets of $\mathcal{D}.$ A proof can be found in Section 11.3 of \cite{B-M}. An application of Gallagher's lemma and Cauchy's integral formula for the derivative of $\zeta$ yields a discrete version of (\ref{m.v.a.c.}), namely 
\begin{align}\label{disc}
\lim\limits_{m\to\infty}\limsup\limits_{N\to\infty}\dfrac{1}{N}\sum\limits_{n=1}^N|\zeta(s+ix_n)-\zeta_m(s+ix_n)|^2=0
\end{align}
uniformly on compact subsets of $\mathcal{D},$ where $(x_n)_{n\in\mathbb{N}}$ is an increasing sequence of non-negative real numbers such that $x_n=O(n)$ and $x_{n+1}-x_n=\Omega(1),$ for all $n\in\mathbb{N}$. 

From the theory of Bergman spaces (see Theorem 1 in \cite{DS}), we know that if $U\subset\mathbb{C}$ is a bounded domain and $f:U\to\mathbb{C}$ is a holomorphic function  which is square integrable on $U$, then
\begin{align}\label{Ber}
|f(z)|\leq\left(\sqrt{\pi}d\left(z,\partial U\right)\right)^{-1}\left(\int_U|f(s)|^2\mathrm{d}A(s)\right)^{1/2}
\end{align}
for all $z\in U$, where $d\left(z,\partial U\right)=\min\left\{|z-w|:w\in\partial U\right\}$ and $\mathrm{d}A$ is the planar Lebesgue measure. Hence, if $K$ is a compact subset of $\mathcal{D}$, we choose a bounded domain $U$ such that $K\subset U\subset\overline{U}\subset\mathcal{D}$ and then it follows from (\ref{disc}) and (\ref{Ber}) that 
$$\lim\limits_{m\to\infty}\limsup\limits_{N\to\infty}\dfrac{1}{N}\sum\limits_{n=1}^N\left(\max\limits_{s\in K}|\zeta(s+ix_n)-\zeta_m(s+ix_n)|\right)^2=0.$$
From the latter equation we can obtain the following lemma.
\noindent
\begin{Lemma}\label{l1}If $(x_n)_{n\in\mathbb{N}}$ is an increasing sequence of non negative real numbers, where $x_n$ satisfy the conditions mentioned above, and $K$ a compact subset of $\mathcal{D},$ then for every $\varepsilon>0$ there exists positive integer $m_0$ such that, for every $m\geq m_0$, 
$$
\liminf\limits_{N\to\infty}\dfrac{1}{N}\#\left\{1\leq n\leq N:\max\limits_{s\in K}|\zeta(s+ix_n)-\zeta_m(s+ix_n)|<\varepsilon\right\}>1-\varepsilon.
$$
\end{Lemma}

\noindent For the needs of this note, we present a simplification of a lemma that is proven in \cite{P}.

\begin{Lemma}\label{l2}Assume that the sequences of real numbers $(x_n)_{n\in\mathbb{N}}$ and $(y_n)_{n\in\mathbb{N}}$ are such that the sequence 
$$
\left(\left(x_n\dfrac{\log p}{2\pi}\right)_{p\in M_1},\left(y_n\dfrac{\log p}{2\pi}\right)_{p \in M_2}\right)\qquad\mbox{for}\quad n\in\mathbb{N},
$$ 
is uniformly distributed mod 1 for any finite sets of primes $M_1$ and $M_2$. Moreover, let $K$ be a compact subset of $\mathcal{D}$ with connected complement and $f,g$ continuous non-vanishing functions on $K,$ which are analytic in the interior of $K.$ Then, for every $\varepsilon>0,$ there exists $m_1>0$ such that for every $m\geq m_1$ we have 
$$
\liminf\limits_{N\to\infty}\dfrac{1}{N}\#\left\{1\leq n\leq N:\begin{array}{ll}\max\limits_{s\in K}|\zeta_{m}(s+ix_n)-f(s)|<\varepsilon\\
\max\limits_{s\in K}|\zeta_{m}(s+iy_n)-g(s)|<\varepsilon \end{array}\right\}>c
$$
with suitable constant $c>0,$ which does not depend on $m_1.$
\end{Lemma}

\noindent For the notion of uniform distribution modulo one we refer to the original work \cite{weyl} of Weyl. 
\par

It is clear that if the sequences $x_n=t_1+\delta_1\lfloor n\alpha\rfloor$ and $y_n=t_2+\delta_2\lfloor n\alpha'\rfloor,$ $n\in\mathbb{N},$ satisfy (Lemma \ref{l1} and) Lemma \ref{l2}, then Theorem \ref{thana} holds. So it suffices to prove that for $\alpha\in\mathcal{L}(\delta_1,\delta_2)\cap(1,+\infty)$, the sequence 
$$
\left(\left(\left(t_1+\lfloor n\alpha\rfloor\delta_1\right)\dfrac{\log p}{2\pi}\right)_{p\in M_1},\left(\left( t_2+\lfloor n\alpha'\rfloor\delta_2\right)\dfrac{\log p}{2\pi}\right)_{p\in M_2}\right)\qquad\mbox{for}\quad n\in\mathbb{N},
$$
is uniformly distributed $\bmod\,1$ for any finite sets of primes $M_1$ and $M_2$. By Weyl's criterion (see Theorem 6.2 of Chapter 1 in \cite{KN}) it is equivalent to show that 
$$
\lim\limits_{N\to\infty}\dfrac{1}{N}\sum\limits_{n=1}^{N}\exp\left[2\pi i\left(\lfloor n\alpha\rfloor\delta_1\sum\limits_{p\in M_1}k_p\dfrac{\log p}{2\pi }+\lfloor n\alpha'\rfloor\delta_2\sum\limits_{p\in M_2}\ell_p\dfrac{\log p}{2\pi }\right)\right]=0,
$$ 
where $(k_p)_{p\in M_1}$ and $(\ell_p)_{p\in M_2}$ are integers and at least one of them is non-zero, or 
$$
\lim\limits_{N\to\infty}\dfrac{1}{N}\sum\limits_{n=1}^{N}\exp\left[2\pi i\left(\lfloor n\alpha\rfloor
\theta_1+\lfloor n\alpha'\rfloor\theta_2\right)\right]=0,
$$ 
where $(\theta_1,\theta_2)\in\mathcal{A}$ and
\begin{eqnarray*}
\lefteqn{\sum\limits_{n=1}^{N}\exp\left[2\pi i\left(\lfloor n\alpha\rfloor\theta_1+\lfloor n\alpha'\rfloor\theta_2\right)\right]}\\
&=&\sum\limits_{n=1}^{N}\exp\left[2\pi i\left(n\left(\alpha\theta_1+\alpha'\theta_2\right)-\lbrace n\alpha\rbrace\theta_1-\lbrace n\alpha'\rbrace \theta_2\right)\right]\\
&=&\sum\limits_{n=1}^{N}f\left(n\left(\alpha\theta_1+\alpha'\theta_2\right),n\alpha,n\alpha'\right),
\end{eqnarray*}
with the Riemann integrable function $f$ defined by 
$$
f(x_1,x_2,x_3)=\exp\left[2\pi i \left(x_1-\lbrace x_2\rbrace\theta_1-\lbrace x_3\rbrace\theta_2\right)\right],
$$
for $(x_1,x_2,x_3)\in\mathbb{R}^3$; notice that $f$ is $1$-periodic in each coordinate and $\lbrace x\rbrace$ denotes the fractional part of the real number $x$. Recall that the sequence 
$$
\mathbf{x}_n=\left(n\left(\alpha\theta_1+\alpha'\theta_2\right),n\alpha,n\alpha'\right)\qquad\mbox{for}\quad 
n\in\mathbb{N},
$$ 
is uniformly distributed mod 1, since $(\theta_1,\theta_2)\in\mathcal{A}$ and $\alpha\in\mathcal{L}(\delta_1,\delta_2)\cap(1,+\infty)$. From Theorem 6.1 of Chapter 1 in \cite{KN}, where the assumption of $f$ being continuous can be substituted with $f$ being Riemann integrable, it follows that
\begin{eqnarray*}
\lefteqn{\lim\limits_{N\to\infty}\dfrac{1}{N}\sum\limits_{n=1}^{N}\exp\left[2\pi i\left(\lfloor n\alpha\rfloor\theta_1+\lfloor n\alpha'\rfloor\theta_2\right)\right]}\\
&=&\lim\limits_{N\to\infty}\dfrac{1}{N}\sum\limits_{n=1}^{N}f(\mathbf{x}_n)=\int\limits_{[0,1]^3}f(\mathbf{x})\mathrm{d}\mathbf{x}\\
&=&\int_0^1\exp(2\pi ix_1)\mathrm{d}x_1\int_0^1\int_0^1\exp[-2\pi i (\theta_1\lbrace x_2\rbrace+\theta_2\lbrace x_3\rbrace)]\mathrm{d}x_2\mathrm{d}x_3
=0.
\end{eqnarray*}
The proof of Theorem \ref{thana} is now almost complete. It only remains to answer the question of existence of such an $\alpha$. For that purpose let $\mathcal{P}\left(\mathbb{R}\right)$ be the power set of $\mathbb{R}$ and define a function 
$$
\ell:\left(\mathbb{Z}^4\setminus\lbrace\textbf{0}\rbrace\right)\times\mathcal{A}\to
\mathcal{P}\left(\mathbb{R}\right)
$$
by mapping the vector $(\textbf{k},\boldsymbol{\theta})=(k_1,k_2,k_3,k_4,\theta_1,\theta_2)$ to the set of $x\in\R$ satisfying
\begin{align}\label{me1}
\left(k_2+k_4\theta_1\right)x^2+\left(k_1-k_2+k_3-k_4\theta_1+k_4\theta_2\right)x-k_1=0.
\end{align}
Observe that, for a given vector in the domain-set of $\ell$, the set described by (\ref{me1}) contains at most two real numbers, which means that it is a set of zero Lebesgue measure contained in $\R$. Thus, 
$$
C=\bigcup\limits_{(\textbf{k},\boldsymbol{\theta})\in\left(\mathbb{Z}^4\setminus\lbrace\textbf{0}\rbrace\right)\times\mathcal{A}}\ell(\textbf{k},\boldsymbol{\theta})
$$
is a subset of $\R$ having zero measure as countable union of sets of zero measure, and $\mathbb{R}\setminus\mathcal{L}(\delta_1,\delta_2)=C$. This proves Theorem \ref{thana}.

Finally, Corollary \ref{sis} follows directly from the just proven theorem by using the fact that, for a given real number $\alpha>1$ and $A\subseteq\mathbb{N}$,  
$$
\liminf\limits_{N\to\infty}\dfrac{1}{N}\#\left\{1\leq n\leq N:n\in A\right\}
\geq \dfrac{1}{\alpha}\liminf\limits_{N\to\infty}\dfrac{1}{N}\#\left\{1\leq n\leq N:\lfloor n\alpha\rfloor\in A\right\}.
$$
Hence, if 
$$
A_\varepsilon=\left\{n\in\mathbb{N}:\begin{array}{ll}\max\limits_{s\in K}|\zeta(s+i(t_1+\delta_1n))-f(s)|<\varepsilon\\
\max\limits_{s\in K}|\zeta(s+i(t_2+\delta_2\sigma_\alpha(n)))-g(s)|<\varepsilon\end{array}\right\},
$$
then 
\begin{eqnarray*}
\lefteqn{\liminf\limits_{N\to\infty}\dfrac{1}{N}\#\left\{1\leq n\leq N:n\in A_\varepsilon\right\}}\\
&\geq&\dfrac{1}{\alpha}\liminf\limits_{N\to\infty}\dfrac{1}{N}\#\left\{1\leq n\leq N:\lfloor n\alpha\rfloor\in A_\varepsilon\right\}\\
&=&\dfrac{1}{\alpha}\liminf\limits_{N\to\infty}\dfrac{1}{N}\#\left\{1\leq n\leq N:\begin{array}{ll}
\max\limits_{s\in K}|\zeta(s+i(t_1+\delta_1\lfloor n\alpha\rfloor))-f(s)|<\varepsilon\\
\max\limits_{s\in K}|\zeta(s+i(t_2+\delta_2\lfloor n\alpha'\rfloor))-g(s)|<\varepsilon
\end{array}\right\},
\end{eqnarray*}
which is positive by Theorem \ref{thana}.

\section{Once More the Right Half of the Critical Strip and Beyond the Critical Line}

Recall that ${\mathcal D}$ denotes the open right half of the critical strip and let $\mathcal{H}(\mathcal{D})$ be the space of holomorphic functions on $\mathcal{D}$ endowed with the topology of uniform convergence on compacta. More precisely, we consider the topology induced from the metric $\rho$ on $\mathcal{H}(\mathcal{D})$ defined by
$$
\rho(f, g) = \sum_{j = 0}^{\infty}  { ||f - g||_{\mathsf{K}_{j}} \over 2^j ( 1 + ||f - g||_{\mathsf{K}_j} ) } 
$$
for any $f, g \in \mathcal{H}(\mathcal{D})$, where $\{ \mathsf{K}_j \}_{j  = 0}^{\infty}$ is a sequence of compact subsets of $\mathcal{D}$ such that $\mathsf{K}_{j} \subset \mathsf{K}_{j + 1}$,
$$
\mathcal{D} = \bigcup_{j = 0}^{\infty} \mathsf{K}_j,
$$ 
and $|| f - g ||_A = \sup_{z \in A} \{ |f(z) - g(z)| \}$. The topology induced from $\rho$ is equivalent to the topology of uniform convergence on compacta. In particular, this implies that the topology of $\mathcal{H}(\mathcal{D})$ is metrizable.

Denote by $\gamma$ the unit circle $\{ z \in \mathbb{C} \,:\, |z| = 1\}$ and put 
$$
\Omega = \prod_{p} \gamma_p
$$
where $\gamma_p = \gamma$ for each prime $p$. With the product topology and operation of point-wise multiplication, this infinite-dimensional torus $\Omega$ is a compact topological Abelian group. Therefore, on $(\Omega, \mathcal{B}(\Omega))$, where $\mathcal{B}(X)$ denotes the Borel $\sigma$-field of $X$, there exists a probability Haar measure $\mathsf{m}_H$. This gives a probability space $(\Omega, \mathcal{B}(\Omega), \mathsf{m}_H)$. Denote by $\omega(p)$ the projection of an element $\omega \in \Omega$ to the coordinate space $\gamma_p$ for each prime $p$ and on the probability space $(\Omega, \mathcal{B}(\Omega), \mathsf{m}_H)$ define the $\mathcal{H}(\mathcal{D})$-valued random element $\zeta(z, \omega)$ by
$$
\zeta(z, \omega) = \prod_p \left( 1 - {\omega(p) \over p^{z} } \right)^{-1}.
$$
Recall that $\zeta(z, \omega)$ converges uniformly on compact subsets of $\mathcal{D}$ for almost all $\omega\in\Omega$ (see \cite[Lemma 4.2]{Ste07}).
Let $\mathsf{P}$ be the distribution of $\zeta(z, \omega)$, that is,
$$
\mathsf{P}(A) = \mathsf{m}_H(\{ \omega \in \Omega \mid \zeta(z, \omega) \in A \}) 
$$
for each $A \in \mathcal{B}(\mathcal{H}(\mathcal{D}))$. We refer to \cite[Chapter $4$]{Ste07} for more details.

A useful tool to prove Theorem \ref{AZTHM1} is the following discrete limit theorem:

\begin{Lemma}\label{DLT}\cite[Proposition $4.4.1$]{Bag81}
For each $h \in \mathbb{R}_{>0}$, the probability measure $\mathsf{P}_N$ defined by
$$
\mathsf{P}_N (A) = {1 \over N} \# \{ n \in \{1, 2, \cdots, N \} \mid \zeta(z + i h n ) \in A \}
$$ 
for each $N \in \mathbb{N}$ and $A \in \mathcal{B}(\mathcal{H}(\mathcal{D}))$ converges weakly to $\mathsf{P}$ as $N \rightarrow \infty$, that is, for any open set $U \subset \mathcal{H}(\mathcal{D})$,
$$
\liminf_{N \rightarrow \infty} \mathsf{P}_N(U) \geq \mathsf{P}(U).
$$
\end{Lemma}

\noindent We remark that Dubickas and Laurin\v cikas \cite[Theorem $7$]{DL16} proved a little different form of Lemma \ref{DLT}: {\it For each $h \in \mathbb{R}_{>0}$ and $\alpha \in (0, 1)$, the probability measure defined by 
$$
\mathsf{P}_{N, \alpha}(A) = {1 \over N } \#\{ n \in \{ 1, 2, \cdots, N \} \mid \zeta(z + i h n^{\alpha}) \in A \}
$$
for each $N \in \mathbb{N}$ and $A \in \mathcal{B}(\mathcal{H}(\mathcal{D}))$ converges weakly to $\mathsf{P}$ as $N \rightarrow \infty$.}
Repeating the same argument, one can replace $n$ with $n^{\alpha}$ and $\alpha$ being any positive quantity in Theorem \ref{AZTHM1} and Corollary \ref{AZTHM2}.

It is also useful to give an explicit form of the support of $\mathsf{P}$, the subset of the domain of $\mathsf{P}$ which contains all points in the domain that do not vanish under $\mathsf{P}$.

\begin{Lemma}\label{ST}\cite[Lemma $5.12$]{Ste07}
The support $\textnormal{Supp}(\mathsf{P})$ of the probability measure $\mathsf{P}$ is 
$$
\{ f \in \mathcal{H}(\mathcal{D}) \mid 1/f \in \mathcal{H}(\mathcal{D}), \; \text{or} \; f \equiv 0\}.
$$
\end{Lemma}

Now we are ready to prove Theorem \ref{AZTHM1} and Corollary \ref{AZTHM2}. We begin with the proof of Theorem \ref{AZTHM1} by using the discrete universality. 

Let $a \in \mathsf{M}$ be an interior point except for $\infty$ and we take $\epsilon>0$ satisfying $\mathsf{B}(a;\epsilon)\subset\mathsf{M}$, where $\mathsf{B}(a;\epsilon):=\{z \in \mathbb{C} \mid |z - a| < r\}$.
Let
$$
\mathsf{K} = \{ s + i h (k - 1) \mid k \in \{ 1, 2, \cdots, l \} \} \subset \mathcal{D}
$$ 
be a finite set with length $l$. Put
$$
U = \{ f \in \mathcal{H}(\mathcal{D}) \mid \sup_{z \in \mathsf{K}} \{ |f(z) - a | \} < \epsilon \}.
$$
Here, we consider the constant $a$ as a constant function defined on $\mathcal{D}$. It is clear that $U \subset \mathcal{H}(\mathcal{D})$ is an open neighborhood of $a \in \mathcal{H}(\mathcal{D})$. It follows from Lemma \ref{ST} that $a \in \textnormal{Supp}(\mathsf{P})$. Moreover, by Lemma \ref{DLT}, we have
$$
\liminf_{N \rightarrow \infty} \mathsf{P}_N(U) \geq \mathsf{P}(U) > 0.
$$
Put $\mathsf{N} = \{ n \in \mathbb{N} \mid \zeta(z + i h n ) \in U \}$. Then $\mathsf{N}$ is infinite and 
$$
|\zeta(s + i h (n + k - 1) )  - a| < \epsilon
$$ 
for any $n \in \mathsf{N}$ and $k \in \{1, 2, \cdots, l  \}$, i.e.,
$$
\zeta(s + i h (n + k - 1) ) \in \mathsf{M}.
$$

Next we prove Corollary \ref{AZTHM2} by using the functional equation for the zeta-function, i.e., 
\begin{equation}\label{feq}
\zeta(s) = \chi(s) \zeta(1- s),
\end{equation}
where
$$ 
\chi(s) = 2^s \pi^{s - 1} \sin{ \left( {\pi s \over 2} \right) } \Gamma(1 - s). 
$$
We begin with an estimate for $\chi$ on the left half of the critical strip.
For convenience, we write $s=\sigma+it$ for $\sigma, t\in\mathbb{R}$.

\begin{Lemma}\label{AZLEM1}
For any $c>0$, there exists a sufficiently large $t_0\geq2$ such that for any
$s=\sigma+it$ with $\sigma\in(0,1/2)$ and $t\geq t_0$, we have
\begin{equation}\label{eq:AZLEM1}
|\chi(s)| \geq c.
\end{equation}
\end{Lemma}

\noindent To see that let $\sigma\in(0,1/2)$ and $t\geq2$.
Applying Stirling's formula to $\Gamma(1-s)$, we have
$$
\log{\Gamma(1-s)}
= \left(\frac{1}{2}-s\right)\log{(1-s)} - 1 + s + \frac{1}{2}\log{2\pi} + O\left(\frac{1}{t}\right)
$$
(see \cite[page 524]{MV06}, for instance).
Thus
$$
\log{|\Gamma(1-s)|} = \operatorname{Re}\left( \log{\Gamma(1-s)} \right)
= \left(\frac{1}{2}-\sigma\right)\log{t} - \frac{\pi t}{2} + \frac{1}{2}\log{2\pi} + O\left(\frac{1}{t}\right).
$$
By expressing sine function in terms of exponential functions, we can easily show that
$$
\log{\left| \sin{\left( {\pi s \over 2} \right)} \right|} = \frac{\pi t}{2} - \log{2} + O(e^{-\pi t})
$$
for $t\geq2$.
We can then show
\begin{align*}
\log{|\chi(s)|} &= \log{| 2^s \pi^{s - 1} \sin{\left( {\pi s \over 2} \right)} \Gamma(1 - s) |} \\
&= \sigma\log{2} + (\sigma-1)\log{\pi} + \frac{\pi t}{2} - \log{2} + O(e^{-\pi t}) \\
&\quad\quad+ \left(\frac{1}{2}-\sigma\right)\log{t} - \frac{\pi t}{2} + \frac{1}{2}\log{2\pi} + O\left(\frac{1}{t}\right) \\
&= (\sigma-1)\log{2\pi} + \left(\frac{1}{2}-\sigma\right)\log{t} + \frac{1}{2}\log{2\pi} + O\left(\frac{1}{t}\right) \\
&= \left(\frac{1}{2}-\sigma\right)\log{\frac{t}{2\pi}} + O\left(\frac{1}{t}\right).
\end{align*}
Hence we can find a sufficiently large $t_0$ such that for any $t\geq t_0$,
$$
\log{|\chi(s)|} > \log{c}
$$
for all $\sigma\in(0,1/2)$. \\

We conclude this section with the proof of Corollary \ref{AZTHM2}.
Let $h>0, l \in \mathbb{N}$, $\mathsf{M} \in \mathbb{M}_{\infty}$ and $\sigma\in(0,1/2)$.
Again denote by $\mathsf{B}(a; r)$ the region $\{z \in \mathbb{C} \mid |z - a| < r\}$.
Without loss of generality, we may assume that $\mathsf{M} \sqcup \mathsf{B}(0; r) = \hat{\mathbb{C}}$ for some $r > 0$.
Applying Theorem \ref{AZTHM1} to $\hat{\mathbb{C}}\backslash\mathsf{B}(0; 2r/c)$, we obtain an infinite subset $\mathsf{N'} \subset \mathbb{N}$ satisfying
$$
\zeta \left( 1 - s - i h (n + k -1 )  \right) \notin \mathsf{B}(0; 2r/c)
$$
for any $n \in \mathsf{N'}$ and $k \in \{1, 2, \cdots, l \}$.
Take a subset $\mathsf{N}$ of $\mathsf{N'}$ such that $n_0 := \min\{n\in\mathsf{N}\}$ satisfies
$t+hn_0 \geq t_0$.
In view of (\ref{feq}) this yields 
\begin{align*}
|\zeta(s + i h (n + k - 1 ))| 
&= |\chi( s + i h (n + k - 1 ) ) | \cdot | \zeta(1 - s - i h (n + k - 1 ) )| \\
&\geq c \cdot 2r/c > r
\end{align*}
for all $n\in\mathsf{N}$ and $k \in \{1, 2, \cdots, l\}$.
Therefore $\zeta(s + i h (n + k - 1 )) \in \mathsf{M}$ for all $n\in\mathsf{N}$ and $k \in \{1, 2, \cdots, l\}$ and this completes the proof.

\section{Concluding Remarks}

It appears that we have no result for the critical line so far. The so-called Hardy $Z$-function is defined as 
$$
Z(t)=\zeta({\textstyle{1\over 2}}+it)\exp(i\theta(t)),
$$
where $\exp(i\theta(t))=\chi({1\over 2}+it)^{-1/2}$ and, consequently,
$$
\theta(t)={t\over 2}\log{t\over 2\pi e}-{7\over 8}+O\left({1\over t}\right),
$$
as $t\to\infty$. It is well-known that $Z(t)$ is real-valued for real $t$ and infinitely often differentiable with the same absolute value as $\zeta({\textstyle{1\over 2}}+it)$. This yields some kind of polar coordinates representation. In fact, if {\bf 耕} denotes $Z(t)$ and {\bf 二} stands for the reciprocal of $\exp(i\theta(t))$, then {\bf 耕二} equals $\zeta({\textstyle{1\over 2}}+it)$. Therefore, we can reformulate the identity from Theorem \ref{absconv}, $\zeta(s+i(t_1+\delta_1n))=\zeta(s+i(t_2+\delta_2\sigma(n)))$ as 
$$
Z(t_1+\delta_1n)=Z(t_2+\delta_2\sigma(n))\qquad \&\qquad \theta(t_1+\delta_1n)\equiv \theta(t_2+\delta_2\sigma(n))\bmod\,2\pi.
$$
Now one may be tempted to use the above approximation for $\theta$ for the congruences $\bmod\,2\pi$. This already leads to strong diophantine conditions for the imaginary parts such that it is rather unlikely that the latter conditions on $Z$ and $\theta$ can hold for many integers $n$. We hope to consider this question more deeply in the near future.  
\bigskip

\noindent {\bf Acknowledgements.} The authors wish to express their gratitude to the organizers of the conference {\it Various Aspects of Multiple Zeta Functions} in honor of Kohji Matsumoto's 60th birthday at Nagoya University in August 2017. Moreover, the authors wish to thank the anonymous referee for his or her careful review.
This work was partly supported by JSPS KAKENHI Grant Numbers 15J02325 and 16J01139.
Further, a part of this work was done under RIKEN Special Postdoctoral Researcher program.

\bigskip

\footnotesize 
\noindent Junghun Lee,
\begin{addmargin}[5mm]{0cm}
\quad Graduate School of Mathematics, Nagoya University, Furo-cho, Chikusa-ku, Nagoya 464-8602, Japan, m12003v@math.nagoya-u.ac.jp; \\
\quad 7th Company, 2nd Battalion, 7th Regiment, 6th Infantry Divisions, 5th Corps, 3rd ROK Army, Republic of Korea
\end{addmargin}
\smallskip

\noindent Athanasios Sourmelidis,
\begin{addmargin}[5mm]{0cm}
\quad Department of Mathematics, W\"urzburg University, Emil-Fischer-Str. 40, 97\,074 W\"urzburg, Germany, athanasiossourmelidis@mathematik.uni-wuerzburg.de
\end{addmargin}
\smallskip

\noindent J\"orn Steuding,
\begin{addmargin}[5mm]{0cm}
\quad Department of Mathematics, W\"urzburg University, Emil-Fischer-Str. 40, 97\,074 W\"urzburg, Germany, steuding@mathematik.uni-wuerzburg.de
\end{addmargin}
\smallskip

\noindent Ade Irma Suriajaya,
\begin{addmargin}[5mm]{0cm}
\quad RIKEN Interdisciplinary Theoretical and Mathematical Sciences Program (iTHEMS) Special Postdoctoral Researcher, RIKEN, 2-1 Hirosawa, Wako, Saitama 351-0198, Japan, adeirmasuriajaya@riken.jp
\end{addmargin}

\end{CJK}
\end{document}